\documentclass[12pt,prl,aps,amssymb,amsmath,tightenlines]{revtex4} 

\usepackage{dsfont}
\newcommand{\half}{\mbox{$\textstyle{\frac{1}{2}}$}}
\newcommand{\re}{{\rm e}}
\newcommand{\ri}{{\rm i}}
\newcommand{\rd}{{\rm d}}

\begin{document}
\title{Operator-valued zeta functions and Fourier analysis}
\author{Carl M. Bender$^1$ and Dorje C. Brody$^{2}$}

\affiliation{
$^1$Department of Physics, Washington University, St Louis, MO 63130, USA\\
$^2$Department of Mathematics, University of Surrey, Guildford GU2 7XH, UK, and
\\
Department of Optical Physics and Modern Natural Science, St Petersburg
National Research University of Information Technologies, Mechanics and Optics, 
St Petersburg 197101, Russia}

\begin{abstract}
The Riemann zeta function $\zeta(s)$ is defined as the infinite sum $\sum_{n=1}^\infty n^{-s}$, which converges when ${\rm Re}\,s>1$. The Riemann hypothesis asserts that the nontrivial zeros of $\zeta(s)$ lie on the line ${\rm Re}\,s= \half$. Thus, to find these zeros it is necessary to perform an analytic continuation to a region of complex $s$ for which the defining sum does not converge. This analytic continuation is ordinarily performed by using a functional equation. In this paper it is argued that one can investigate some properties of the Riemann zeta function in the region ${\rm Re}\,s<1$ by allowing {\it operator-valued} zeta functions to act on test functions. As an illustration, it is shown that the locations of the trivial zeros can be determined purely from a Fourier series, without relying on an explicit analytic continuation of the functional equation satisfied by $\zeta(s)$.
\end{abstract}

\maketitle

\noindent {\bf 1}. 
In the so-called Hilbert-P\'olya programme, one attempts to establish the
Riemann hypothesis by (a) finding an operator (possibly a differential operator)
along with a boundary condition such that the eigenvalues of the operator
correspond to the nontrivial zeros of the zeta function, and (b) showing that
the operator is self-adjoint. A ${\cal PT}$-symmetric operator fulfilling the
requirement (a) has recently been identified \cite{r1} and its eigenvalues were
subsequently shown to be real \cite{r2}. However, these findings have not yet
brought us closer to a proof of the Riemann hypothesis because the span of the
eigenfunctions of the operator satisfying the boundary condition may not be
fully contained within the rigged Hilbert space upon which the self-adjointness
of the operator is established. To be specific, although a Hilbert space has
been identified, one is automatically restricting the analysis to the critical
line ${\rm Re}\,s=\half$ of the zeta function $\zeta(s)$ so that little can be
inferred about the zeros off the critical line, if there are any. This may
explain why the Hilbert-P\'olya programme has not yet led to a proof of the
Riemann hypothesis.

With this in mind, we propose here an alternative way to investigate the
properties of the Riemann zeta function by analysing {\it operator-valued} zeta
functions; that is, we examine the behaviour of $\zeta(\hat s)$, where $\hat s$
is an operator, such as the dilation operator. It is unclear whether this
approach will lead to a deep understanding of the zeta function. However, it is
shown here that one can use this approach to establish some properties of
Riemann zeta function very easily. Specifically, we can calculate $\zeta(s)$
for some integer values of $s$ without performing an explicit analytic
continuation.

In the Hilbert-P\'olya programme one investigates operators whose eigenvalues 
correspond to the locations of the zeros of the zeta function but we propose
here to investigate operators whose eigenvalues are the values of the zeta
function itself. For example, in the Hilbert-P\'olya programme one might
consider the properties of the operator $\half(1-\ri{\hat h}_{\rm BK})$, where
${\hat h}_{\rm BK}={\hat x}{\hat p}+{\hat p}{\hat x}$ denotes the Berry-Keating
Hamiltonian \cite{r3}, but here we consider the operator $\zeta(\half(1-\ri{\hat
h}_{\rm BK}))$. We investigate such an operator by letting it act on
trigonometric functions. We will show that
\begin{equation}
\zeta\left(\half(1-\ri{\hat h}_{\rm BK})\right)\sin x=\frac{\sin x}{2(1-\cos x)}
\label{eq:1}
\end{equation} 
for $x\in(0,\pi]$, from which we can deduce that $\zeta(-n)=(-1)^nB_{n+1}/(n+1)$
for $n$ a positive-odd integer, where $\{B_n\}$ are the Bernoulli numbers. There
are numerous similar relations, from which further properties of the zeta
function can be inferred. Another example is 
\begin{equation}
\zeta\left(\half(3-\ri{\hat h}_{\rm BK})\right)\sin x=\frac{\pi-x}{2} 
\label{eq:2}
\end{equation} 
for $x\in[0,\pi]$, from which we can deduce that $\zeta(s)$ vanishes for
negative-even integers $s$ (these are the trivial zeros) without analytically
continuing the functional equation $\zeta(s)=2^s\pi^{s-1}\sin(\pi s/2)\Gamma(1-s
)\zeta(1-s)$. From (\ref{eq:2}) we can also deduce that $\zeta(0)=-\frac{1}{2}$,
and that $\zeta(s)$ has a pole at $s=1$. 
We remark that operators of the form 
$\zeta(\half\pm\half\ri{\hat h}_{\rm BK})$ were mentioned briefly in Ref.~\cite{r3}.  

\vspace{0.2cm}
\noindent {\bf 2}. 
This paper is based in part on the following relations in Fourier analysis
\cite{r4}: 
\begin{equation}
\sum_{n=1}^\infty\frac{\cos(nx)}{n^{2m}}=\frac{(-1)^{m-1}(2\pi)^{2m}}{2(2m)!}\,
B_{2m}\left(\frac{x}{2\pi}\right)  
\label{eq:3}
\end{equation} 
and 
\begin{equation}
\sum_{n=1}^\infty\frac{\sin(nx)}{n^{2m-1}}=\frac{(-1)^{m}(2\pi)^{2m-1}}
{2(2m-1)!}\,B_{2m-1}\left(\frac{x}{2\pi}\right)  
\label{eq:4}
\end{equation} 
for $m=1,2,\ldots$, 
where $B_m(x)$ denotes the Bernoulli polynomial of order $m$. Note that these
equations are only valid for {\it real} $x\in[0,2\pi]$; the Fourier series
diverge for complex $x$. If we set $m=1$ in (\ref{eq:4}), we obtain $(\pi-x)/2$,
which is the right side of (\ref{eq:2}). Similar series were investigated by
Clausen \cite{r5}. 

In this paper we reinterpret these results by using quantum-mechanical operator
techniques. We show that it is possible to infer properties of $\zeta(s)$ by
studying actions of the operator $\zeta({\hat s})$ on functions when ${\hat s}$
is an operator. (Bernoulli polynomials play a major role in the theory of the
Riemann zeta function, so it is not surprising that some properties of the zeta
function on the real line can be inferred in the context of Fourier series.) 

Series of the form (\ref{eq:3}) or (\ref{eq:4}) can be extended to cases for
which $m$ is $0$ or a negative integer. Such series can be summed by using Euler
summation. For example, we have 
\begin{equation}
\sum_{n=1}^\infty\re^{{\rm i}nx}=\lim_{r\to1^-}\sum_{n=1}^\infty\left(r\re^{{\rm
i}x}\right)^n=\lim_{r\to1^-}\frac{r\re^{{\rm i}x}}{1-r\re^{{\rm i}x}}=
\frac{1}{\re^{-{\rm i}x}-1}. 
\label{eq:5}
\end{equation} 
Taking the imaginary part, we deduce that 
\begin{equation}
\sum_{n=1}^\infty\sin(nx)=\frac{\sin x}{2(1-\cos x)}, 
\label{eq:6}
\end{equation} 
which is the right side of (\ref{eq:1}). This result can also be obtained by 
complex analysis; from the analytic continuation of the Lerch zeta function 
\begin{equation}
L(s,x)=\sum_{n=1}^\infty\frac{\re^{{\rm i}nx}}{n^s}=\frac{\Gamma(1-s)}{2\pi\ri}
\int_C\frac{\re^t\,t^{s-1}}{1-\re^{t+{\rm i}x}}\,\rd t,
\label{eq:7}
\end{equation}
Apostol deduced sums such as (\ref{eq:6}) \cite{r6}. Here, the integration path
$C$ is a Hankel contour that encircles the negative-$t$ axis in the positive
direction. Thus, we obtain (\ref{eq:5}) by setting $s=0$ in (\ref{eq:7}) and
using the residue at $t=0$ to evaluate the integral. 

\vspace{0.2cm}
\noindent {\bf 3}. 
Let us proceed to establish relations such as (\ref{eq:2}) above. For this
purpose we require the notion of the dilation operator. The generator of the
dilation is ${\hat x}{\hat p}$, where ${\hat p}=-\ri\,\rd/\rd x$, so that for a
smooth function $f(x)$ we have 
$$\re^{{\rm i}\lambda{\hat x}{\hat p}}f(x)=f(\re^\lambda x).$$
It follows that 
$$\sin(nx)=n^{{\rm i}{\hat x}{\hat p}}\,\sin x.$$
Therefore, ignoring for now the question of the convergence of the sum, we
deduce that
\begin{equation}
\sum_{n=1}^\infty\frac{\sin(nx)}{n}=\sum_{n=1}^\infty\frac{n^{{\rm i}{\hat x}{
\hat p}}}{n}\,\sin x=\zeta\left(1-\ri{\hat x}{\hat p}\right)\,\sin x.
\label{eq:10}
\end{equation} 
Thus, the action of the Riemann dilation operator $\zeta(1-\ri{\hat x}{\hat p})$
on a trigonometric function generates a Fourier series. 
[Note that we do not define an operator of the form $\zeta(z+\ri{\hat x}{\hat p})$ 
as a Taylor expansion of $\zeta(s)$ about $s=z$ in powers of 
$\ri{\hat x}{\hat p}$. Such an expansion may diverge.] 
In this example,
the left side is the Fourier representation for the linear function $(\pi-x)/2$,
and from the relation ${\hat h}_{\rm BK}=2{\hat x}{\hat p}-\ri$ we observe that
$1-\ri{\hat x}{\hat p}=\half(3-\ri{\hat h}_{\rm BK})$. We therefore deduce the
identity (\ref{eq:2}). Hence if the operator $\zeta\left(1-\ri{\hat x}{\hat p}
\right)$ were invertible, we would expect the relation 
$$\frac{1}{\zeta\left(1-\ri{\hat x}{\hat p}\right)}\,\frac{\pi-x}{2}=\sin x$$
to hold. However, since ${\hat x}{\hat p}$ is the dilation generator, it cannot
change the power of $x$ on the left side, so we arrive at a contradiction. This 
suggests that the operator $\zeta\left(1-\ri{\hat x}{\hat p}\right)$ {\it cannot be 
inverted} because its spectrum contains at least one zero eigenvalue. 

Before we proceed to inspect the locations of the zeros, let us check the 
consistency of (\ref{eq:2}) without relying on the summation representation of
the zeta function. For this purpose we use the integral representation 
\begin{equation} 
\zeta(s)=\frac{\Gamma(1-s)}{2\pi\ri}\int_C\frac{t^{s-1}}{\re^{-t}-1}\,\rd t 
\label{eq:12}
\end{equation}
for the zeta function and 
$$\frac{1}{\Gamma(1-s)}=\frac{1}{2\pi\ri}\int_C\re^{t}t^{s-1}\,\rd t$$
for the reciprocal of the Gamma function. Because $s$ appears in two different 
ways in (\ref{eq:12}) our strategy is to check the validity of 
$$\frac{1}{\Gamma(\ri{\hat x}{\hat p})}\,\frac{\pi-x}{2}=\frac{\zeta\left(1-\ri{
\hat x}{\hat p}\right)}{\Gamma(\ri{\hat x}{\hat p})}\,\sin x $$
to infer (\ref{eq:2}). For the left side we deduce that 
$$\frac{1}{\Gamma(\ri{\hat x}{\hat p})}\,\frac{\pi-x}{2}=\frac{1}{2\pi\ri}
\int_C\re^{t}t^{-{\rm i}{\hat x}{\hat p}}\left(\frac{\pi-x}{2}\right)\rd t=\frac
{1}{2\pi\ri}\int_C\re^t\left(\frac{\pi-t^{-1}x}{2}\right)\rd t=-\frac{1}{2}x.$$
The constant term $\pi/2$ has been annihilated here because of the pole of 
$\zeta(s)$ at $s=1$. On the other hand, expanding $\sin x$ in a power series, 
we deduce from 
$$\frac{\zeta\left(1-\ri{\hat x}{\hat p}\right)}{\Gamma(\ri{\hat x}{\hat p})}\, 
x^n=\frac{1}{2\pi\ri}\int_C\frac{t^{-{\rm i}{\hat x}{\hat p}}}{\re^{-t}-1}x^n\,
\rd t=\frac{1}{2\pi\ri}\int_C \frac{t^{-n}}{\re^{-t}-1}x^n\,\rd t
=\frac{\zeta(1-n)}{\Gamma(n)}\,x^n$$
that 
\begin{equation}
\frac{\zeta\left(1-\ri{\hat x}{\hat p}\right)}{\Gamma(\ri{\hat x}{\hat p})}\, 
\sin x=\sum_{n=1}^\infty\frac{\zeta(2(1-n))}{(2n-1)!\Gamma(2n-1)}\, x^{2n-1}.   
\label{eq:17}
\end{equation}
Since the right side of (\ref{eq:17}) must equal $-\frac{1}{2}x$, we infer that 
$\zeta(0)=-\frac{1}{2}$, and that $\zeta(-2)=\zeta(-4)=\cdots=0$. Conversely,
from these elementary facts about the zeta function we infer the consistency of
(\ref{eq:10}). 

An essentially identical line of argument leads to the observation that 
\begin{equation}
\zeta\left(2-\ri{\hat x}{\hat p}\right)\,\cos x=
\frac{\pi^2}{6}-\frac{\pi x}{2}+\frac{x^2}{4},
\label{eq:18}
\end{equation}
\begin{equation}
\zeta\left(3-\ri{\hat x}{\hat p}\right)\,\sin x=
\frac{\pi^2x}{6}-\frac{\pi x^2}{4}+\frac{x^3}{12},
\label{eq:19}
\end{equation}
and so on. Thus, for each of the Clausen functions in (\ref{eq:3}) and
(\ref{eq:4}) we obtain a corresponding representation in the form of an operator
$\zeta\left( N-\ri{\hat x}{\hat p}\right)$ acting on a trigonometric function,
for $N$ a positive integer. Each of these relations reveals some information
about the values of $\zeta(s)$ for real integral values of $s$.

\vspace{0.2cm}
\noindent {\bf 4}. 
As a slightly shorter way to do the analysis above, we observe
that since $\ri{\hat x}{\hat p}\,x^\alpha=\alpha x^\alpha$, and since $\zeta(s)$
is analytic except for a simple pole at $s=1$, we have $\zeta(N-\ri{\hat x}{\hat
p})\,x^n=\zeta(N-n)\,x^n$. However, one must be careful about the existence of
the pole. To illustrate this, we consider the example $\zeta(1-\ri{\hat x}{\hat
p})\,\sin x$. Expanding the sine series, and assuming the interchangeability of
the two limits, we obtain 
$$\zeta(1-\ri{\hat x}{\hat p})\sum_{n=1}^\infty\frac{x^{2n-1}}{(2n-1)!}=\sum_{
n=1}^\infty\zeta(1-\ri{\hat x}{\hat p})\frac{x^{2n-1}}{(2n-1)!}=\sum_{n=1
}^\infty\zeta(2-2n)\frac{x^{2n-1}}{(2n-1)!}=-\frac{1}{2}x,$$
which shows that term-by-term application of the differential operator $\zeta(
1-\ri{\hat x}{\hat p})$ is not permissible because we have missed the constant
term $\pi/2$ associated with the pole of $\zeta(s)$. In fact, for each of the
examples discussed above, interchanging the limits leaves out just one term
corresponding to the pole of $\zeta(s)$; that is, one term on the left side that
is annihilated by $\Gamma(\ri{\hat x}{\hat p}+1-N)^{-1}$. This term is the only
parity-violating term; while each term on the left side of (\ref{eq:4}) has odd
parity, one term on the right side has even parity. Similarly, while each term
on the left side of (\ref{eq:3}) has even parity, one term on the right side has
odd parity. Thus, in (\ref{eq:18}) the term $\pi x/2$ on the right side violates
parity, and similarly in (\ref{eq:19}) the term $\pi x^2/4$ on the right side
violates parity. Hence, the commutator of the two limits gives the
parity-breaking term resulting from summing the series. 

We remark that the term-by-term application of the operator respects both parity
and analyticity. To illustrate this, we consider the series in (\ref{eq:6}).
Observe that each term in the series on the left side has odd parity and the
right side is also odd, so that there is no violation of parity. One might
expect that term-by-term application of $\zeta(-\ri{\hat x}{\hat p})$ on the
power-series expansion of $\sin x$ is permissible. However, while each term in
the series on the left side is analytic and vanishes at $x=0$, the right side 
diverges like $1/x$ as $x\to0$. Indeed,
\begin{eqnarray}
\sum_{n=0}^\infty\zeta(-\ri{\hat x}{\hat p})\frac{(-1)^n}{(2n+1)!}x^{2n+1} &=&
\sum_{n=0}^\infty\zeta(-2n-1)\frac{(-1)^n}{(2n+1)!}\,x^{2n+1}\nonumber\\
&=& \sum_{n=0}^\infty(-1)^{2n+1}\frac{B_{2n+2}}{2n+2}\frac{(-1)^n}{(2n+1)!}\,
x^{2n+1}\nonumber\\ 
&=& \frac{1}{x}\sum_{n=0}^\infty\frac{\ri^{2n+2}}{(2n+2)!}\,B_{2n+2}\,x^{2n+2} 
\nonumber\\
&=&\frac{1}{x}\sum_{k=2}^\infty\frac{1}{k!}\,B_k\,(\ri x)^{k}\nonumber\\
&=& \frac{1}{x}\left[\sum_{k=0}^\infty\frac{1}{k!}\,B_k\,(\ri x)^k-1-\ri B_1x
\right].
\label{eq:21}
\end{eqnarray}
Therefore, from the generating function $\sum_{k=0}^\infty B_k\,x^k/k!=x/(\re^x-
1)$ with $B_1=-\frac{1}{2}$ we get
$$\sum_{n=0}^\infty\zeta(-\ri{\hat x}{\hat p})\frac{(-1)^n}{(2n+1)!}x^{2n+1}=
\frac{\sin x}{2(1-\cos x)}-\frac{1}{x}.$$
Remarkably, we recover the right side of (\ref{eq:1}), but with its singularity 
removed. Moreover, the singular term in the right side of (\ref{eq:1})
corresponds to the pole of $\zeta(s)$ at $s=1$. This is the only term that is
annihilated by the action of $\Gamma(1+\ri{\hat x}{\hat p})^{-1}$.

Analogous results can be seen in other examples, for instance, in
$$\zeta\left(-1-\ri{\hat x}{\hat p}\right)\,\cos x=\sum_{n=1}^\infty n\cos(nx)=-
\frac{1}{2(1-\cos x)}.$$
Once again, there is no parity violation but the right side is singular at $x=0$
and behaves like $-1/x^2$, while each of the summands in the middle term is well
behaved. On the other hand, by interchanging the order of differentiation and
summation associated with the Taylor expansion of $\cos x$ we obtain 
$$\sum_{n=0}^\infty\zeta\left(-1-\ri{\hat x}{\hat p}\right)\frac{(-1)^n}{(2n)!} 
x^{2n}=-\sum_{n=0}^\infty\frac{(2n+1)(\ri x)^{2n}}{(2n+2)!}\,B_{2n+2}.$$
Then a calculation like that in (\ref{eq:21}) leads to the same conclusion that
$$\sum_{n=0}^\infty\zeta\left(-1-\ri{\hat x}{\hat p}\right)\frac{(-1)^n}{(2n)!} 
x^{2n}=-\frac{1}{2(1-\cos x)}+\frac{1}{x^2},$$
and the singularity at the origin has been removed. 

\vspace{0.2cm}
\noindent {\bf 5}. 
The analysis presented here can be extended to more general Dirichlet 
$L$-functions. These are functions expressible in the form 
$$
L_\chi(s) = \sum_{n=1}^\infty \frac{\chi(n)}{n^s} 
$$ 
for ${\rm Re}(s)>1$, and otherwise can be defined by their analytic continuations. 
Here $\chi(n)$ denotes a Dirichlet character, which is a function from integers to 
complex numbers satisfying the multiplicative property that $\chi(mn)=\chi(m)
\chi(n)$, the periodicity that $\chi(n)=\chi(n+k)$ for some positive $k$, and the 
condition that if $n$ and $k$ are relative primes then $\chi(n)\neq0$ but otherwise 
$\chi(n)=0$. Thus, for $k=1$ we have $\chi(n)=1$ for all $n$ and $L_\chi(s)$ 
reduces to the Riemann zeta function. 

As a simple example other than the Riemann zeta function, let us consider the 
Dirichlet beta function arising from considering the period $k=4$. Specifically, for 
${\rm Re}(s)>1$ the Dirichlet beta function is defined by the series 
$$ 
\beta(s) = \sum_{n=0}^\infty \frac{(-1)^n}{(2n+1)^s} ,
$$ 
from which we deduce that 
$$ 
\beta(-\ri{\hat x}{\hat p})\, \sin x  = \sum_{n=0}^\infty (-1)^n 
(2n+1)^{{\rm i}{\hat x}{\hat p}} \sin x 
= \sum_{n=0}^\infty (-1)^n \sin\big( (2n+1)x\big)
=0,
$$ 
where the vanishing of the alternating sine series here can be deduced by using 
Euler summation. On the other hand, interchanging the order of summation and 
differentiation in the series expansion of $\sin x$ gives 
$$ 
\beta(-\ri{\hat x}{\hat p})\, \sin x  = \sum_{n=0}^\infty \frac{(-1)^n}{(2n+1)!} \, 
\beta\big(-(2n+1)\big) \, x^{2n+1}, 
$$ 
from which we deduce that $\beta(-n)=0$ for all positive odd $n$, without explicitly 
relying on analytic continuation. Note that the 
interchange of the limits is permissible in this example because there is no pole 
contribution. 

An analogous calculation shows that 
$$ 
\beta(-\ri{\hat x}{\hat p})\, \cos x  = \sum_{n=0}^\infty (-1)^n 
(2n+1)^{{\rm i}{\hat x}{\hat p}} \cos x 
= \sum_{n=0}^\infty (-1)^n \cos\big( (2n+1)x\big)
=\frac{1}{2\cos x},
$$ 
whereas by interchanging the limits we find that 
$$ 
\beta(-\ri{\hat x}{\hat p})\, \cos x  = \sum_{n=0}^\infty \frac{(-1)^n}{(2n)!} \, 
\beta(-2n) \, x^{2n}. 
$$ 
Comparing these two we deduce that $\beta(-n)=E_n/2$ for all positive even 
$n$. This result can also be obtained by considering 
$$ 
\beta(1-\ri{\hat x}{\hat p})\, \sin x  = \sum_{n=0}^\infty \frac{(-1)^n}{2n+1}  
\sin\big( (2n+1)x\big)=\half\ri\big[\tan^{-1}(\re^{-{\rm i}x})-\tan^{-1}(\re^{{\rm i}x})
\big] , 
$$ 
and comparing this with 
$$ 
\beta(1-\ri{\hat x}{\hat p})\, \sin x  = \sum_{n=0}^\infty \frac{(-1)^n}{(2n+1)!} \, 
\beta(-2n) \, x^{2n+1} . 
$$

\vspace{0.2cm}
\noindent {\bf 6}. 
In the foregoing analysis we have only considered one class of operator-valued 
zeta functions, namely, zeta functions evaluated at a linear function of the 
dilation operator. This class of operators is suitable in the context of Fourier
analysis \cite{r7}. It appears that the action of this class of operators on
trigonometric functions only yields information about $\zeta(s)$ for real $s$
although further study is required to clarify this point. In this connection, we
note that the matrix elements of, for example, $\zeta(1-\ri{\hat x}{\hat p})$,
viewed as an operator acting on the Hilbert space of square-integrable functions
on $[0,\pi]$, in the standard sine basis $\{\sqrt{2/\pi}\sin(nx)\}$, is given by
$$\zeta_{mn}=\left\{\begin{array}{ll} n/m & \quad{\rm if}~n~{\rm divides}~m,\\
0 & \quad{\rm otherwise}.\end{array}\right.$$
Thus, the matrix $\{\zeta_{mn}\}$ encodes the information about factorisation of
integers. This suggests that it might be possible to extract more information
by studying further properties of the class of operator-valued zeta functions
considered here. 

Evidently, there are many other operator-valued zeta functions that one might
consider. For instance, the action of $\zeta({\hat p}^2+{\hat x}^2)$ on Hermite
polynomials might yield further results on the zeta function. As another
example, if we let ${\hat a}=({\hat x}+\ri{\hat p})/\sqrt{2}$ denote the
standard annihilation operator and $|s\rangle$, $s\in{\mathds C}$, a coherent
state, we then have $\zeta({\hat a})|s\rangle=\zeta(s)|s\rangle$. Thus, if the
action of the operator $\zeta({\hat a})$ were implementable in a laboratory,
then one would see the coherent light being absorbed whenever $s$ is a zero of
the zeta function. 

To conclude, we have shown that by studying the action of Riemann dilation
operators on trigonometric functions, we are able to infer some properties of
the Riemann zeta function. Of course, the properties of $\zeta(s)$ inferred here
are already known. Nevertheless, we were able to determine, for example, the
locations of the trivial zeros from elementary Fourier analysis without relying
explicitly on the analytic continuation of the zeta function. This suggests that
further research into actions of operator-valued zeta functions may yield
interesting new results.

\vspace{0.2cm} 
\noindent 
{\bf Acknowledgement}

\vspace{0.1cm} 
\noindent 
DCB thanks the Russian Science Foundation for support (project 16-11-10218). The 
authors thank J.~Keating for suggesting the idea of examining other Dirichlet 
$L$-functions.

\end{document}